\documentclass[12pt,leqno]{article}

\newcounter{conjecture}
\newcounter{remark}
\newcounter{corollary}

\newenvironment{corollary}{\medskip{\bf Corollary \thecorollary.}
\addtocounter{corollary}{1}\em}{\rm}
\newtheorem{theorem}{Theorem}
\newtheorem{lemma}{Lemma}

\newcommand {\rrr}[1]{(\ref{#1})}

\def \be{\begin{equation}}
\def \ee{\end{equation}}
\def \bt{\begin{theorem}}
\def \et{\end{theorem}}
\def \bc{\begin{corollary}}
\def \ec{\end{corollary}}
\def \bea{\begin{eqnarray}}
\def \eea{\end{eqnarray}}
\def \bas{\begin{eqnarray*}}
\def \eas{\end{eqnarray*}}
\def \bl{\begin{lemma}}
\def \el{\end{lemma}}

\def \vski{\vspace{12pt}}

\def \({\left(}
\def \){\right)}

\def \bc{\begin{center} }
\def \ec{\end{center} }
\def \bs{\begin{slide} }
\def \es{\end{slide} }

\def\square{{\vcenter{\vbox{\hrule height.3pt
         \hbox{\vrule width.3pt height5pt \kern5pt
            \vrule width.3pt}
         \hrule height.3pt}}}}

\newcounter{cccases}
\begin{document}

\title{An electric-resistance approach to return time}

\author{
\begin{tabular}{c}
\textit{Greg Markowsky} \\
gmarkowsky@gmail.com \\
+61 03 9905-4487 \\
Monash University \\
Department of Mathematical Sciences \\
Victoria, 3800 Australia
\end{tabular}}



\maketitle

\begin{abstract}
A new proof is given for the formula for the expected return time of a random walk on a graph. This proof makes use of known relationships between electric resistance and random walks.
\end{abstract}

\vski


\vski


\vski






Let $G$ be a finite, connected graph with $n$ vertices and $m$ edges. We will use the standard notation $\sim$ to denote adjacency in the graph, and for $z \in G$ let $deg(z)$ denote the degree of $z$. Let $X_j$ denote simple random walk on $G$; that is, $X_j$ is the Markov chain taking values in the vertex set of $G$ with transition probabilities given by

\be P(X_{j+1} = z | X_j = y) = \left \{ \begin{array}{ll}
\frac{1}{deg(y)} & \qquad  \mbox{if } y \sim z  \\
0 & \qquad \mbox{otherwise } \;.
\end{array} \right. \ee

Let $T_z = \inf\{j \geq 0 : X_j =z \}$, and let $T_z^+ = \inf\{j \geq 1 : X_j =z \}$. The object of interest for us is the expected return time, $E_z[T^+_z]$. The following elegant theorem is well known.

\begin{theorem} \label{t1}
\begin{equation} \label{}
E_z[T^+_z] = \frac{2m}{deg(z)}.
\end{equation}
\end{theorem}

This result admits a considerable generalization. Consider each edge $(y,z)$ as a wire in a circuit with a given conductance $C_{yz}$, which is a nonnegative number which measures how easily electricity (and the random walk) passes along the edge. For each vertex $z$ let $C_z = \sum_{y \sim z} C_{yz}$. Let $X_j$ now be the Markov chain taking values in the vertex set of $G$ with transition probabilities given by

\be P(X_{j+1} = z | X_j = y) = \left \{ \begin{array}{ll}
\frac{C_{yz}}{C_y} & \qquad  \mbox{if } y \sim z  \\
0 & \qquad \mbox{otherwise } \;.
\end{array} \right. \ee

This is the random walk induced by the electric network, and simple random walk corresponds to taking conductances of 1 (or any positive constant) across each edge. It should be noted that this construction is in fact quite general, since any reversible Markov chain can be realized as such an induced random walk (see \cite[Ch. 9]{levin}). Let $C = \sum_{y \in G} C_y$. We then have the following extension of Theorem \ref{t1} (which in fact applies to infinite graphs as well under the assumption that $C$ is finite).

\begin{theorem} \label{t2}

\begin{equation} \label{sin}
E_z[T^+_z] = \frac{C}{C_z}.
\end{equation}
\end{theorem}

The standard method of proving Theorem \ref{t2} is to appeal to a result from Markov chain theory, namely that an irreducible Markov chain with a stationary distribution $\pi$ satisfies $E_z[T^+_z] = 1/\pi_z$; and then simply verifying that $\pi_z = \frac{C_z}{C}$ is the stationary distribution for the chain $X_j$ (see \cite[Sec. 1.7]{norris}). On the other hand, researchers studying electric resistance have uncovered many identities and bounds on such quantities as hitting times, commute times, and cover times (\cite{comcov}); mixing times (\cite[Ch. 4]{aldfill}); and edge-cover times (\cite{george}). It is therefore natural to search for a derivation of Theorems \ref{t1} and \ref{t2} which makes more use of the principles which relate electric resistance to random walks, especially in light of the statement of the second theorem. We now present such a proof, naturally of the more general result, Theorem \ref{t2}.

\vski

Fix $z \in G$, and construct a new graph $\tilde G$ which contains $G$ as a subgraph by adding a vertex $\tilde z$ to $G$ with a single edge connecting $\tilde z$ to $z$. Across the new edge lay a conductance of 1. Let $\tilde X_m$ denote a random walk on $\tilde G$ induced by the conductances present (the original ones in $G$, together with the edge with unit conductance connecting $z$ and $\tilde z$). For $x,y \in \tilde G$, let $\tilde R_{x,y}$ be the effective resistance within $\tilde G$ between $x$ and $y$, and let $\tilde T_y = \inf\{j \geq 0 : \tilde X_j =y \}$. Let $\tilde C$ be twice the sum of the conductances across all edges in $\tilde G$; note that $\tilde C=C+2$. It is known (\cite[Cor. 11, Ch. 3]{aldfill}) that

\begin{equation} \label{}
E_{\tilde z}[\tilde T_z] + E_z[\tilde T_{\tilde z}] = \tilde C \tilde R_{z,\tilde z}.
\end{equation}

However, it is trivial that $E_{\tilde z}[\tilde T_z] = 1$ (the walk beginning at $\tilde z$ has no choice but to pass to $z$ at time 1), and it is equally trivial that $R_{z,\tilde z} = 1$ (there are no paths from $\tilde z$ to $z$ except for along the edge connecting them). Making the necessary substitutions yields

\begin{equation} \label{hai}
E_z[\tilde T_{\tilde z}] = C+1.
\end{equation}

Now, at time $\tilde T_{\tilde z} - 1$ the walk $\tilde X_m$ must necessarily reside at $z$. Furthermore, at each visit to $z$ the walk $\tilde X_m$ has probability $\frac{1}{C_z+1}$ of passing to $\tilde z$ at the next step, and between visits to $z$ the walk performs excursions within $G$, which will each take an average of $E_z[T^+_z]$ steps. The number of excursions within $G$ before $\tilde T_{\tilde z}$ is a Bernoulli trial with probability $\frac{1}{C_z+1}$, and as is well known the expected number of such trials until first success is $C_z+1$; however the number of excursions will in fact be one less than the number of visits to $z$, since the walk begins at $z$. It follows then that the expected number of excursions within $G$ before $\tilde T_{\tilde z}$ will be $C_z$. Adding $1$ to record the final step from $z$ to $\tilde z$ we obtain

\begin{equation} \label{shang}
E_z[\tilde T_{\tilde z}] = C_z E_z[T^+_z]+1.
\end{equation}

Equating the right-hand sides of \rrr{shang} and \rrr{hai} yields \rrr{sin}.

\vski

{\bf Remark:} The key idea of attaching a new vertex to a vertex $z$ in a graph and starting a random walk there, armed with the knowledge that the first step must be to $z$, appears also in different contexts in \cite[Lemma 3.1]{northpal} and \cite[Ex. 10.4]{levin}.

\section{Acknowledgements}

I'd like to thank Tim Garoni and Jos\'e Palacios for helpful conversations. I am also grateful for support from Australian Research Council Grants DP0988483 and DE140101201.

\bibliographystyle{alpha}
\bibliography{prob}

\end{document}